\documentclass[aap,MSNbibl,seceqn,citesort,dvips]{arximspdf}
\usepackage{graphics}
\doi{10.1214/10-AAP689}
\volume{21}
\issue{1}
\pubyear{2011}
\firstpage{115}
\lastpage{139}

\makeatletter
\newtheorem{theorem}{Theorem}[section]
\newtheorem{lemma}{Lemma}[section]
\def\re{{\mathbb{R}}}
\def\Rset{{\mathbb{R}}}

\def\e{{\mathbb{E}}}
\def\ve{{\varepsilon}}
\def\T{\Phi}
\def\H{\mathcal{H}}
\def\R{\mathbf{r}}
\makeatother

\begin{document}
\begin{frontmatter}

\title{A general framework for waves in random media with long-range
correlations}
\runtitle{Waves in long-range random media}

\begin{aug}
\author[a]{\fnms{Renaud} \snm{Marty}\ead[label=e1]{renaud.marty@iecn.u-nancy.fr}\corref{}}
\and
\author[b]{\fnms{Knut} \snm{S{\o}lna}\ead[label=e2]{ksolna@math.uci.edu}}
\address[a]{Institut Elie Cartan de Nancy\\
Nancy-Universit\'e, CNRS,
INRIA\\
B.P. 239, F-54506 Vandoeuvre-l\`es-Nancy Cedex\\
France\\
\printead{e1}}
\address[b]{University of California at Irvine\\
Irvine, California 92697-3875\\
USA\\
\printead{e2}}
\affiliation{Universit\'e Henri Poincar\'e Nancy 1 and University of
California at Irvine}

\runauthor{R. Marty and K. S{\o}lna}
\end{aug}

% HISTORY:
\received{\smonth{6} \syear{2009}}
\revised{\smonth{2} \syear{2010}}

% ABSTRACT
%
\begin{abstract}
We consider waves propagating in a randomly layered medium with
long-range correlations.
An example of such a medium is studied in \cite{MS} and leads, in
particular, to an asymptotic
travel time described in terms of a fractional Brownian motion. Here we
study the asymptotic
transmitted pulse \mbox{under} very general assumptions on the
long-range correlations. In the framework
that we introduce in this paper, we prove in particular that the
asymptotic time-shift can be described in terms of non-Gaussian
and/or multifractal processes.
\end{abstract}

% KEYWORDS
%


\end{frontmatter}

%s1 ###
\section{Introduction}
Wave propagation in random media has been extensively studied for many
years from both
theoretical and applied points of view.
In particular, the study of the effective shape of an acoustic pulse
propagating through a layered medium has attracted a lot of attention
\cite{akppw,CF,SolPap}.
Recently, applications to time reversal \cite{book1} have also been the
subject of much interest.
Currently, there is also a strong interest in problems related to noise
and correlations \cite{garnier05}. In all
these situations the statistical properties of the medium are
important since they affect the statistical properties of the wave field.

In \cite{CF} the authors consider an acoustic pulse propagating in
a\break
one-dimensional random medium with rapidly decaying correlations. They
rigorously prove
the classical result of O'Doherty and Anstey \cite{oda} which
establishes that
the effective transmitted pulse is characterised by deterministic
spreading and a random time-shift. More precisely, the deterministic spreading
is expressed as a convolution with a Gaussian density, and the
random time-shift is described in terms of a Brownian motion.

More recently, wave propagation and also homogenization
in random media with long-range correlations and/or
defined in terms of fractional Brownian motion \cite{BGMP,GS,MS,Sol}
have been considered.
In \cite{MS}, we extend the result of \cite{CF}
to such a framework. Then, the asymptotic description of the transmitted
pulse is dramatically different from what happens in a mixing case. Indeed,
the pulse keeps its initial shape, and its random time-shift is now
described in terms of a fractional Brownian motion whose Hurst index depends
on the decay rate of the correlation function of the random fluctuations.
We considered in \cite{MS} a particular form of a random
process describing the medium, such that it was roughly speaking close
to a Gaussian process.
Thus, it still remains to study more general cases under the long-range
assumption.
This is the aim of the present work. We establish that under general
long-range assumptions on the medium, the effective pulse still keeps
its initial
shape as observed in \cite{MS}, but the time-shift can be very
different, non-Gaussian,
for instance, depending on the form of the random fluctuations.
Moreover, our general
result allows us to deal with media with
a decay rate of correlations that varies along the propagation direction.
This leads to an
effective time-shift, described in terms of a multifractional random
process which is,
roughly speaking, a fractional Brownian motion with a varying Hurst
index that reflects
the nonhomogeneity of the propagation medium.

In Section \ref{sec1} we introduce the problem and review
the basic wave decomposition approach.
Next, we establish the general technical result (Theorem \ref{main}) in
Section~\ref{sec2}.
We apply
this general result
to non-Gaussian media
in Section \ref{nonGaussian},
to multifractal Gaussian media in Section \ref{multifractal}
and to multifractal
non-Gaussian media in Section \ref{nonGm}
where we prove the main result of the paper (Theorem \ref{main4}).
We present a numerical illustration in Section \ref{secnum}.
Finally, Section \ref{secmain} is devoted to the derivation
of Theorem~\ref{main}.

%s2 ###
\section{Preliminaries}\label{sec1}
%s2.1 ###
\subsection{Wave decomposition}
The governing equations are the nondimensionalized Euler equations
giving conservation of moments
and mass
%
%e2 ###
%e1 ###
\begin{eqnarray}
\rho^{\ve} (z){\partial u^{\ve}\over\partial t}(z,t)+{\partial p^{\ve
}\over\partial z} (z,t)& = & 0  , \\
{1\over K^{\ve}(z)}{\partial p^{\ve}\over\partial t}(z,t)+{\partial
u^{\ve}\over\partial z} (z,t)& = & 0  ,
\end{eqnarray}
where $t$ is the time, $z$ is the depth into the medium, $p^{\ve}$ is the
pressure and $u^{\ve}$ the particle velocity. The medium parameters are
the density
$\rho^{\ve}$ and the bulk-modulus $K^{\ve}$ (reciprocal of the compressibility).
We assume that $\rho^{\ve} $ is a constant identically equal to one
in our nondimensionalized setting,
and $1/K^{\ve}$ is modeled as follows:
%
%e3 ###
\begin{equation}
{1\over K^{\ve}(z)}=\cases{
{1+\mu^{\ve}({z})}
&\quad  $\mbox{for } z\in[0,Z],$\vspace*{2pt}\cr
1 & \quad $\mbox{for }z\in\re- [0,Z],$}
\end{equation}
where $\mu^{\ve}$ is a centered random process.
The number $\ve>0$ is a parameter that all quantities depend on. As we
will see below it is introduced to describe the scales of the problem.

We introduce the right- and left-going waves
%
%e4 ###
\begin{equation}
A^{\ve}=p^{\ve}+u^{\ve} \quad \mbox{and} \quad
B^{\ve}=u^{\ve}-p^{\ve}.
\end{equation}
The boundary conditions are of the form
%
%e5 ###
\begin{equation}
A^{\ve}(z=0,t) = f(t/\ve^{\tau})\quad  \mbox{and}\quad  B^{\ve}(z=Z,t)=0
\end{equation}
for a positive real number $\tau>0$ and a source function $f$.
In order to deduce a description of the transmitted pulse, we
open a window of size $\ve^{\tau}$ in the neighborhood of the
travel time of the homogenized medium
and define the processes
%
%e6 ###
\begin{equation}\label{decomp}
a^{\ve}(z,s) = A^{\ve}(z,z+\ve^{\tau}s) \quad  \mbox{and}\quad  b^{\ve}(z,s) =
B^{\ve}(z,-z+\ve^{\tau}s) .
\end{equation}
Observe that the background or homogenized medium
in our scaling has a constant speed of sound equal to unity
and that the medium is matched so that in the frame introduced
in (\ref{decomp}) the pulse shape is constant if $\mu^{\ve} \equiv0$
or if we consider the homogenized medium \cite{book1}.
We introduce next
the Fourier transforms $\hat{a}^{\ve}$ and $\hat{b}^{\ve}$ of
$a^{\ve}$
and $b^{\ve}$, respectively,
\[
\hat{a}^{\ve}(z,\omega)=\int_{-\infty}^{\infty} e^{i\omega s}a^{\ve
}(z,s)\,ds \quad  \mbox{and} \quad  \hat{b}^{\ve
}(z,\omega)=\int_{-\infty}^{\infty} e^{i\omega s}b^{\ve}(z,s)\,ds   ,
\]
that satisfy
%
%e8 ###
%e7 ###
\begin{eqnarray}
{d\hat{a}^{\ve}\over dz} &= &{i\omega\over2}
\nu^{\ve} ({z} )
(\hat{a}^{\ve}-e^{-2i\omega z/\ve^{\tau}}\hat{b}^{\ve}
),\qquad\hat{a}^{\ve}(0,\omega)=\hat{f}(\omega)   ,
\\
{d\hat{b}^{\ve}\over dz} &=& {i\omega\over2}
\nu^{\ve}({z})(e^{2i\omega z/\ve^{\tau}}\hat{a}^{\ve
}-\hat{b}^{\ve}),\qquad\hat{b}^{\ve}(Z,\omega)=0,
\end{eqnarray}
where we use the notation
%
%e9 ###
\begin{equation}
\label{nudefinition}
\nu^{\ve} = {\mu^{\ve}\over\ve^{\tau}}.
\end{equation}
Following \cite{CF,book1} we express the previous system of equations
in terms
of the propagator $P^{\ve}_{\omega}(z)$ which can be written as
%
%e10 ###
\begin{equation}
P^{\ve}_{\omega}(z)=
\pmatrix{
\alpha^{\ve}_{\omega}(z) &\overline{\beta^{\ve}_{\omega}}(z)\vspace*{2pt} \cr
\beta^{\ve}_{\omega}(z) & \overline{\alpha^{\ve}_{\omega}}(z)},
\end{equation}
and that satisfies
%
%e11 ###
\begin{equation}
{dP_{\omega}^{\ve}\over dz} (z)= \H_{\omega}^{\ve}
\biggl({z\over\ve^{\tau}},{z}\biggr)P_{\omega}^{\ve}(z),\qquad
P_{\omega}^{\ve}(z=0)=\pmatrix{
1 & 0 \cr
0 & 1} ,
\end{equation}
with
\[
\H_{\omega}^{\ve}(z_1,z_2)={i\omega\over2}\nu^{\ve}(z_2)\pmatrix{
1 & -e^{-2i\omega z_1}\vspace*{2pt}\cr
e^{2i\omega z_1} & -1}.
\]
Defining next the transmission coefficient $T_{\omega}^{\ve}$ and
the reflection coefficient $R_{\omega}^{\ve}$ by
%
%e12 ###
\begin{equation}
T_{\omega}^{\ve}(z)={1\over\overline{\alpha^{\ve}_{\omega}}(z) }
\quad  \mbox{and}\quad
R_{\omega}^{\ve}(z)={\beta^{\ve}_{\omega}(z) \over\overline{\alpha^{\ve
}_{\omega}}(z) },
\end{equation}
we can write
%
%e13 ###
\begin{equation}
\label{representation-a-fourier1}
a^{\ve}(Z,s) =
{1\over2\pi}\int_{-\infty}^{\infty} e^{-is\omega}T_{\omega}^{\ve
}(Z)\hat{f}(\omega)  \,d\omega
\end{equation}
and
%
%e14 ###
\begin{equation}
\label{representation-a-fourier2}
b^{\ve}(0,s) ={1\over2\pi}\int_{-\infty}^{\infty} e^{-is\omega
}R_{\omega}^{\ve}(Z)\hat{f}(\omega)
\,d\omega .
\end{equation}
Hence, we shall study the asymptotics of the propagator $P^{\ve}_{\omega}$
in order to characterize~$a^{\ve}$ and $b^{\ve}$ as $\ve$ goes to 0.
%s2.2 ###
\subsection{A short-range medium}\label{secshort}
We recall now what happens in a
mixing (or short-range) model when $\tau=1$ and $\mu^{\ve}(z)=\nu(z/\ve
^2)$. We assume that
$\nu= \T\circ m$ where $\T$ is a bounded function and $m$ is a
centered Markov process with an invariant probability measure whose
generator satisfies the Fredholm alternative.
This implies that the covariance function $z\mapsto\e[\nu(0)\nu(z)]$
is integrable and then the correlation length $\sigma$ of the medium
is finite
\[
\sigma^2=\int_{0}^{\infty}\e[\nu(0)\nu(z)] \, dz\in[0, \infty)   .
\]
This property is the mixing property or the short-range property.
It is well known
\cite{CF,book1} that under these assumptions the propagator equations
$P_{\omega}^{\ve}$ converge to a system of stochastic differential
equations driven by
independent Brownian motions from which we can deduce that
$a^{\ve}(Z,s) \longrightarrow \widetilde{a}(Z,s) $
as $\ve$ goes to 0 with
%
%e15 ###
\begin{equation}
\widetilde{a}(Z,s)=(f * G)(s-B)   ,
\end{equation}
where $G$ is a centered Gaussian density with variance $\sigma^2 Z/2 $
and $B$ a Gaussian random variable that can be expressed in terms of a
Brownian motion $W$ as $B= \sigma W(Z) / \sqrt{2}$.
Proving this result involves
using the Diffusion Approximation Theorem \cite{book1}
to get an asymptotic propagator from which we can deduce the
expression of the
limit $\widetilde{a}(Z,s)$. Notice that, whereas the variance of $B$
depends in particular on $\T$,
the result does not depend qualitatively on $\T$ in the sense that $B$
remains Gaussian whatever
$\T$ is.
%s2.3 ###
\subsection{A long-range medium}
\label{longrangemedium}
In \cite{MS}, the propagation in a long-range medium is investigated. The
model considered is defined in terms of a fractional Brownian motion.
More precisely, we assume that
$\nu^{\ve}$ has the form
\[
\nu^{\ve}(z)=\ve^{\kappa-\tau}\nu\biggl( {z\over\ve^2}\biggr)\qquad \mbox{for } z\in
[0,Z],
\]
where $\kappa>0$ and $\nu$ is a process that is
expressed as
$\nu(z)=\T(m(z))$
for every $z$
where:
\begin{itemize}
\item$\T$ is an odd $\mathcal{C}^{\infty}$-function;
\item$m$ is a Gaussian process,
centered, stationary and has a correlation function
$r_m$ which has the following asymptotic property as
$z$ goes to $\infty$:
%
%e16 ###
\begin{equation}
\label{assumption}
r_m(z)=\e[m(0)m(z)]\sim c_{m}z^{-\gamma},\qquad
\gamma\in(0,1)   .
\end{equation}
\end{itemize}
The property (\ref{assumption}) implies
that the covariance function $r_{\nu}$ of $\nu$ is not
integrable
\[
\int_{0}^{\infty}|r_{\nu}(z)| \, dz=\infty  ,
\]
which means that the correlation length is infinite.
This is the so-called long-range property.
We mention that a typical example of a process satisfying (\ref
{assumption}) can be
constructed as
%
%e17 ###
\begin{equation}\label{mod1}
m(z)=W_H(z+1)-W_H(z)   ,
\end{equation}
where $B_H$ is a fractional Brownian motion (fBm in short) with Hurst parameter
$H>1/2$.

We assume $\tau$, $\kappa$ and $\gamma$ satisfy $\tau-\kappa= \gamma
$. In this case, we proved that
$a^{\ve}(Z,s) \to \widetilde{a}(Z,s) $
with
%
%e18 ###
\begin{equation}
\label{pulseshapeMS}
\widetilde{a}(Z,s)=f(s-B)   ,
\end{equation}
where $B$ a Gaussian random variable. We can write $B$ as
$B= {\sigma_H}W_H(Z)$ where~$W_H$ is a fractional Brownian motion with
Hurst parameter
$H=(2-\gamma)/2$ and $\sigma_H$ is a positive constant that depends on
$H$ and $\T$.

%s3 ###
\section{Medium assumptions and main technical result}\label{sec2}
The results presented above show that the asymptotic behavior of the
pulse shape $a^{\ve}(Z,s)$ strongly depends on the statistical
properties of $\nu$. The pulse shape is affected under short-range
assumptions whereas it does not change under the long-range assumptions
described above.
In Sections \ref{nonGaussian} and \ref{multifractal} we
carry out the analysis of the particular long-range media
that we consider in this paper.
To facilitate this analysis
we establish in this section a theorem under the following general
assumptions on $\nu^{\ve}=\mu^{\ve}/\ve^{\tau}$.
Let $\lambda>0$ and define:
\begin{itemize}
\item Assumption ${\mathrm{A}_1}$: As $\ve$ goes to 0, the
finite-dimensional distributions of the process
$\{ \int_0^z\nu^{\ve}(z') \, dz' \}_z$ converge to those of a process
$V=\{ V(z)\}_z$ with finite
second-order moments.
\item Assumption ${\mathrm{A}_2(\lambda)}$: There exist two symmetric,
continuous and two-variable functions
$\gamma\dvtx  [0,Z]^2\to[\gamma_-,\gamma_+]\subset(0,1)$ and $R \dvtx
[0,Z]^2\to(0,\infty)$ such that
for every $\delta>0$, there exists $z_{\delta}>0$ sufficiently large
such that
for every $z_1$, $z_2$ and $\ve$ satisfying
$|z_1-z_2|>\ve^{\lambda}z_{\delta}$,
\[
\bigl|\e[\nu^{\ve}(z_1)\nu^{\ve}(z_2)]-R(z_1,z_2)
|z_1-z_2|^{-\gamma(z_1,z_2)}
\bigr|\leq\delta R(z_1,z_2)|z_1-z_2|^{-\gamma(z_1,z_2)}.
\]

\item Assumption ${\mathrm{A}_3(\lambda)}$. For every $\rho>0$ there
exist $C_{\rho}>0$ and $\gamma_{\rho}\in(0,1)$ such that
$|\e[\nu^{\ve}(z_1)\nu^{\ve}(z_2)]|\leq C_{\rho}|z_1-z_2|^{-\gamma_{\rho
}}$ for every $\ve>0$
and $|z_1-z_2| < \ve^{\lambda} \rho$.
\end{itemize}
Assumption ${\mathrm{A}_1}$ corresponds to
the convergence of the travel-times. Assumptions~${\mathrm{A}_2}(\lambda
)$ and ${\mathrm{A}_3}(\lambda)$ are long-range assumptions for
nonstationary processes. They describe how the long-range property
varies with the propagation distance. In particular, these enable
us to apply the next theorem to
multifractal media (Sections~\ref{multifractal} and \ref{nonGm}), which
are nonhomogeneous.

Here we give the main technical result of this paper.
This theorem is next used in Sections \ref{nonGaussian}, \ref{multifractal} and
\ref{nonGm}
 to establish the asymptotic pulse shape respectively in
non-Gaussian and multifractal media.
\begin{theorem}
\label{main}
Assume that there exists $\lambda>0$ such that
${\mathrm{A}_1}$, ${\mathrm{A}_2(\lambda)}$ and ${\mathrm{A}_3(\lambda)}$ are
satisfied.
Then, as $\ve$ goes to 0, $\{a^{\ve}(Z,s)\}_{s}$ converges
in distribution in the space of continuous functions endowed with the
uniform topology to the random process $\{\widetilde{a}(Z,s)\}_{s}$
that can be written as
%
%e19 ###
\begin{eqnarray}
\widetilde{a}(Z,s) = f\bigl(s-\tfrac{1}{2}V(Z) \bigr)   .
\end{eqnarray}
\end{theorem}

Theorem \ref{main} establishes that, under general long-range
assumptions, if the travel-time converges then the asymptotic pulse
keeps its initial shape but its time shift is described in terms of the
asymptotic travel-time.
As recalled in Section \ref{longrangemedium} this fact was observed
in a
particular case in \cite{MS}. In fact, the result of \cite{MS}
follows from Theorem \ref{main}.
Indeed, the model presented in Section \ref{longrangemedium}
satisfies ${\mathrm{A}_1}$, ${\mathrm{A}_2}(2)$ and ${\mathrm{A}_3}(2)$. In
particular,
the finite-dimensional distributions of
$\{ \int_0^z\nu^{\ve}(z')\,  dz' \}_z$ converge to those of the process
$\{ 2\sigma_HW_H(z)\}_z$, so that the asymptotic pulse is of the form
(\ref{pulseshapeMS}).

Notice that the framework we study in this paper is in dramatic
contrast with the mixing case where we observe a pulse spreading in
addition to the time-shift. This is not so
surprising if we remark that Theorem \ref{main} does not apply to a
process~$\nu^{\ve}$ defined as in Section \ref{secshort}
by $\nu^{\ve}(z)=\ve^{-1}\nu(z/\ve^2)$ where $\nu$ is a mixing
process. Indeed, if such a process $\nu^{\ve}$ satisfied assumption
${\mathrm{A}_2(\lambda)}$ for some $\lambda>0$, then we would have
\[
\int_{z^*}^{\infty}\e[\nu^{\ve}(z)\nu^{\ve}(0)]\,dz\geq c^*\int
_{z^*}^{\infty}{dz\over z^{\gamma^*}}=\infty
\]
for some $c^*>0$, $z^*\in[0,Z]$ and $\gamma^*\in(0,1)$, which
contradicts the mixing assumption that gives
\[
\int_{z^*}^{\infty}\e[\nu^{\ve}(z)\nu^{\ve}(0)]\,dz\leq\int_{0}^{\infty
}\e[\nu(z)\nu(0)]\,dz< \infty.
\]

To conclude this section we present a
heuristic description of the link between the mixing and the long-range
cases. For every $\omega$ we define
\[
\mathbf{v}^{\ve}=\mathbf{v}_{\omega}^{\ve}:=(v_1^{\ve},v_{2,\omega}^{\ve
},v_{3,\omega}^{\ve}),
\]
where for every $z\in[0,Z]$ by
\begin{eqnarray*}
v_1^{\ve}(z)&=&\int_0^z
\nu^{\ve}({z'})  \,dz'  ,\\
 v_{2,\omega}^{\ve}(z)&=&\int_0^z
\nu^{\ve}({z'})
\cos\biggl(2\omega{z'\over\ve^{\tau}}\biggr)  \,dz'  ,\\
v_{3,\omega}^{\ve}(z)&=&\int_0^z
\nu^{\ve}({z'})\sin\biggl(2\omega{z'\over\ve^{\tau}}
\biggr)
\,dz'   .
\end{eqnarray*}
In both cases the three-dimensional process $\mathbf{v}^{\ve}$
plays a crucial role. In the mixing case $\mathbf{v}^{\ve}$ converges
to the
three-dimensional (nonstandard) Brownian motion $(B_1,B_{2,\omega
},B_{3,\omega})$. In the proof of the convergence
\[
a^{\ve}(Z,s) \longrightarrow \widetilde{a}(Z,s) =(f * G)(s-B)
\]
one then observes that the Gaussian variable $B$ can be written as $B=B_1(Z)/2$,
and that the Gaussian density $G$ derives from $B_{2,\omega}$ and
$B_{3,\omega}$ \cite{CF,GP,book1}.
In the long-range case, let us assume that $\mathbf{v}^{\ve}$ converges
to the
three-dimensional process $(V,0,0)$. This fact was already observed in
\cite{M2} for the fractional white noise.
Now if we substitute $(B_1,B_{2,\omega},B_{3,\omega})$ with $(V,0,0)$
in the expression of the limit $\widetilde{a}(Z,s) $ we obtain
$B=V(Z)/2$, $G=\delta_0$ (because in fact $\hat{G}\equiv1$) and
hence $\widetilde{a}(Z,s) =f(s-V(Z)/2)$. This is what we establish in
this paper, in particular by proving
the convergence of $\mathbf{v}^{\ve}$ and the substitution mentioned
just above.

%s4 ###
\section{Non-Gaussian asymptotics}
\label{nonGaussian}
In this section we study the case where $\nu^{\ve}$ has the form
\[
\nu^{\ve}(z)=\ve^{\kappa-\tau}\nu\biggl( {z\over\ve^2}\biggr)\qquad  \mbox{for } z\in
[0,Z],
\]
where $\kappa>0$ and $\nu$ is a process that is assumed to
have the form
\[
\nu(z)=\T(m(z))
\]
for every $z$
where:
\begin{itemize}
\item$\T$ is a continuous function such that $\T(\sigma_0 \times\cdot
)$ has a
Hermite index equal to $K\in\mathbb{N}^*$, where $\sigma_0^2=\e[m(0)^2]$.
\item$m$ is a continuous Gaussian process,
centered, stationary and has a correlation function
$r_m$ which has the following asymptotic property as
$z$ goes to $\infty$:
%
%e20 ###
\begin{equation}
\label{assumptionbis}
r_m(z)=\e[m(0)m(z)]\sim c_{m}z^{-\gamma},
\end{equation}
where $0<\gamma<1/K$.
\end{itemize}
We denote the $K$th Hermite coefficient of $\T(\sigma_0 \times\cdot)$ by
\[
J(K)=\e[\T(\sigma_0X)P_K(X)],
\]
where
$X\sim\mathcal{N}(0,1)$, and $P_K$ is the $K$th Hermite polynomial.
Applying Theorem~\ref{main} we get the following result.
\begin{theorem}\label{main2}
Assume that $\tau-\kappa= \gamma K$.
Then, as $\ve$ goes to 0, $\{a^{\ve}(Z,s)\}_{s}$ converges
in distribution in the space of continuous functions endowed with the
uniform topology to the random process $\{\widetilde{a}(Z,s)\}_{s}$
that can be written as
%
%e21 ###
\begin{equation}
\widetilde{a}(Z,s) = f\bigl(s-\tfrac{1}{2}W_H^K(Z) \bigr)   ,
\end{equation}
where $W_H^K$ is the $K$th Hermite process of index $H=(2-\gamma
K)/2\in(1/2, 1)$ defined
for every $z$ by
%
%e22 ###
\begin{equation}
\label{definitionhermite}
W_H^K(z)  =  \frac{c{}^{K/2}_m}{\sigma^{K}_0}
\int_{ \re^K}\mathcal{G}_{H,K}(z,x_1,\ldots, x_K ) \prod
_{k=1}^K\hat{B}(dx_k)
\end{equation}
with
\[
\mathcal{G}_{H,K}(z,x_1,\ldots, x_K )=\frac{}{}
{ J(K) (e^{-iz\sum_{j=1}^Kx_j}-1 )\over K!C(H)^K \sum_{j=1}^Kx_j }
\prod_{k=1}^K{x_k\over|x_k|^{{(H-1)/
K}+3/2}},
\]
where $\hat{B}(dx)$ is the Fourier transform of a Brownian measure,
\[
C(H)^2=\int_{-\infty}^{\infty}\frac{e^{-ix}}{|x|^{1+2(H-1)/K}}\,dx,
\]
and the multiple stochastic integral is in the sense of \cite{Dobrushin}.
\end{theorem}

For $H\in(1/2, 1)$ and $K\in\mathbb{N}^*$ given, the Hermite process
defined by
(\ref{definitionhermite}) was studied independently in \cite{DM} and
\cite{T3}. Its increments
are stationary and its covariance is
\[
\e[W_H^K(z_1)W_H^K(z_2)]=\tfrac{1}{2}(|z_1|^{2H}+|z_2|^{2H}-|z_1-z_2|^{2H}).
\]
It is self-similar and $H$-H\"older.
It is Gaussian if and only if $K=1$; thus, it is a fractional Brownian motion
if and only if $K=1$. As a consequence, the result of~\cite{MS}
corresponds to the case of $K=1$
in Theorem~\ref{main2}. Moreover, this result is in dramatic contrast
to the short-range case
where the asymptotics does not depend qualitatively on $\T$.
\begin{pf*}{Proof of Theorem \protect\ref{main2}}
Following \cite{DM} or \cite{T3}, we find that the finite-dimensional
distributions of the antiderivative
of $\nu^{\ve}$ converge to those of $W_H^K$; therefore,
${\mathrm{A}}_1$ is satisfied.
Next we show that ${\mathrm{A}}_2(2)$ and ${\mathrm{A}}_3(2)$ hold.
Because of the stationarity of $m$ it is enough to show that
%
%e23 ###
\begin{equation}
\label{asymptoticv}
\e[\nu(0)\nu(z)]\sim c_{\nu}z^{-K\gamma}\qquad \mbox{as }z\to
\infty
\end{equation}
for some constant $c_{\nu}>0$.
By the Hermite expansion we can write
\[
\nu(z)=\T\biggl(\sigma_0{m(z)\over\sigma_0}\biggr)=\sum_{k=K}^{\infty
}{J(k)\over k!}P_k\biggl({m(z)\over\sigma_0}\biggr)
.
\]
Using the properties of the Hermite polynomials we get
%
%e24 ###
\begin{eqnarray}\label{asymptoticv2}
\nonumber\e[\nu(0)\nu(z)] & = & \sum_{k=K}^{\infty}{J(k)^2\over
(k!)^2}\e\biggl[
P_k\biggl({m(0)\over\sigma_0}\biggr)P_k\biggl({m(z)\over\sigma_0}
\biggr)\biggr]
\nonumber
\\[-8pt]
\\[-8pt]
\nonumber
& = & \sum_{k=K}^{\infty}{J(k)^2\over k!\sigma
_0^{2k}}r_m(z)^k  .
\end{eqnarray}
Therefore, we need to study the limit of
\[
z^{\gamma K}\e[\nu(0)\nu(z)] = \sum_{k=K}^{\infty}{J(k)^2\over
k!\sigma_0^{2k}}z^{\gamma K}r_m(z)^k  .
\]
Observe that for $k=K$ we have
$z^{\gamma K}r_m(z)\sim c$ as $z \to\infty$, and for $k>K$ we have
$z^{\gamma K}r_m(z)^k\to0$.
Moreover, we have the uniform upper bound for $z$ sufficiently large
\[
{J(k)^2\over k!\sigma_0^{2k}}z^{\gamma K}|r_m(z)|^k\leq{J(k)^2\over
k!}   .
\]
Using the fact that
$
\sum_{k=1}^{\infty} {J(k)^2\over k!}<\infty,
$
(\ref{asymptoticv}) follows from the uniform convergence theorem.
\end{pf*}

%s5 ###
\section{Application to multifractal media}
\label{multifractal}
In this section we study the case where the asymptotic medium is
described in terms of a multifractional
process. In all the situations described above, the media were
asymptotically expressed in terms of fractional
processes. A drawback of fractional processes for
applications is the strong homogeneity of their properties, which are
described by their (constant)
Hurst index. Therefore, multifractional processes have attracted
much attention \cite{BJR,PLV}. Multifractional
processes have locally the same properties as fractional
processes. Their properties are governed
by a $(0,1)$-valued function
$h$ which is called the multifractional function. Some of the main
properties are that
multifractional processes are locally
self-similar, and their pointwise H\"older
exponents vary along their trajectory. In particular, multifractional
processes are relevant in order
to describe nonhomogeneous media. Before stating the main result of
this section, we mention that
the most famous multifractional process is the multifractional Brownian
motion. It was
independently introduced in \cite{BJR,PLV} and can be defined from the
harmonizable representation of fractional Brownian motion for every $z$
%
%e25 ###
\begin{equation}
\label{WH}
W_H(z)=\frac{1}{C(H)}\int_{-\infty}^{\infty}{e^{-izx}-1\over
|x|^{H+1/2}}\hat{B}(dx),
\end{equation}
where $\hat{B}$ is the Fourier transform of a real Gaussian
measure $B$, and the constant $C(H)$ is a renormalization constant and
can be written as
%
%e26 ###
\begin{equation}
\label{CH}
C(H)^2=\int_{-\infty}^{\infty}{|e^{-ix}-1|^2\over|x|^{2H+1}}\,  dx
={\pi\over H\Gamma(2H)\sin(\pi H)}.
\end{equation}
Now we consider a $(0,1)$-valued function $h$, and we substitute $H$ by
$h(z)$ for every~$z$ to obtain
%
%e27 ###
\begin{equation}
\label{Wh}
W_h(z)=\frac{1}{\widetilde{C}(z)}\int_{-\infty}^{\infty
}{e^{-izx}-1\over|x|^{h(z)+1/2}}\hat{B}(dx),
\end{equation}
where the constant $\widetilde{C}(z)$ is a renormalization
function.

We shall here use a different framework for the multifractal
modeling that is convenient for the asymptotic analysis
and describe this next.
We assume that $\nu^{\ve}$ has the form
\[
\nu^{\ve}(z)=\ve^{\kappa(z)-\tau}\nu\biggl( {z\over\ve^2}, z\biggr)\qquad
\mbox{for } z\in
[0,Z],
\]
where $\kappa$ is a positive function, and $\nu$ is a field that is
written as
$\nu(z_1,z_2)=\T(m(z_1,h(z_2)))$
for every $z_1$ and $z_2$
where:
\begin{itemize}
\item$\T$ is a continuous function with Hermite index $1$.
\item$h$ is a continuous function taking values in $[h_-, h_+]\subset
(1/2, 1)$.
\item$m=\{m(z,H)\}_{z, H}$ is a centered and continuous Gaussian field
such that $\e[m(z,H)^2]=1$ for every $z$ and $H$ and
such that
there exists a continuous function $\R\dvtx  [h_-,h_+]^2\to(0,\infty)$
(that we call the asymptotic covariance of $m$) such that
\begin{eqnarray*}
&&\lim_{z_1-z_2\to\infty}\sup_{(H_1, H_2)}| \R( H_1, H_2)\\
&&\hspace*{62pt}\quad {}-(z_1-z_2)^{2-H_1-H_2}\e[m(z_1,H_1)m(z_2,H_2)]|=0.
\end{eqnarray*}
\end{itemize}
These assumptions describe that the field $m$ has the long-range
property with respect to the variable $z$. They also express that for
each $H$, the process $m(\cdot, H)$ is stationary and asymptotically
fractional because it satisfies the
classical invariance principle. As established in \cite{CM} this field
enables us
to define a process that is asymptotically multifractional.

Applying Theorem \ref{main} we now get the following theorem.
\begin{theorem}
\label{main3}
Let $\gamma(z):=\tau-\kappa(z)$, and assume $h(z)=(2-\gamma(z))/2$.
Then, as $\ve$ goes to 0, $\{a^{\ve}(Z,s)\}_{s}$ converges
in distribution in the space of continuous functions endowed with the
uniform topology to the random process $\{\widetilde{a}(Z,s)\}_{s}$
that can be written as
%
%e28 ###
\begin{eqnarray}
\widetilde{a}(Z,s) = f\bigl(s-\tfrac{1}{2}S_h(Z) \bigr)   ,
\end{eqnarray}
where
${S}_h$ is a centered Gaussian process with covariance for $z_1,z_2\geq
0$ given by
%
%e29 ###
\begin{equation}
\label{covlim}
\e[{S}_h(z_1){S}_h(z_2)]  =  J(1)^2
\int_0^{z_1}du_1\int_0^{z_2}du_2\,
\widetilde{\mathcal{R}}( u_1,u_2 ) ,
\end{equation}
where
\[
\widetilde{\mathcal{R}}(u_1,u_2)=\mathcal{R}
(u_1,u_2;h(u_1),h(u_2))|u_1-u_2|^{h(u_1)+h(u_2)-2}
\]
with
%
%e30 ###
\begin{equation}
\label{Rcal}
\mathcal{R}(z_1,z_2;H_1,H_2)=\R( H_1, H_2)1_{z_1\geq z_2}+\R( H_2,
H_1)1_{z_1<z_2}.
\end{equation}
\end{theorem}

The process $S_h$ was introduced in \cite{CM}. This process is
continuous and multifractional in
the sense that its pointwize H\"older exponent is $h(t_0)$ at the point $t_0$:
\[
\sup\biggl\{ H , \lim_{\ve\to0}{S_h(t_0+\ve)-S_h(t_0)\over|\ve|^{H}}=0
\biggr\}=h(t_0).
\]
Notice that in the case of $h$ is constant Theorem \ref{main3}
corresponds to the result of~\cite{MS}.
\begin{pf*}{Proof of Theorem \protect\ref{main3}}
By the same procedure as in proving (\ref{asymptoticv}), we get from the
asymptotic assumptions for
$\{m(z,H)\}$ that
\begin{eqnarray*}
&&\lim_{z_1-z_2\to\infty}\sup_{(H_1, H_2)\in [h_-,h_+]^2}|
(z_1-z_2)^{2-H_1-H_2}\e[\nu(z_1,H_1)\nu(z_2,H_2)]\\
&&\hspace*{184pt}\qquad{}-J(1)^2\R( H_1, H_2) |=0.
\end{eqnarray*}
If we denote, respectively, $v^{\ve}$ and $w^{\ve}$ the antiderivatives of
$z\mapsto\nu^{\ve}(z)$ and $z\mapsto\ve^{2h(z)-2}m(z/ \ve^2,
h(z))$,
then, by using
the same argument as above we also get
%
%e31 ###
\begin{equation}
\lim_{\ve\to0}\e[|v^{\ve}(z)-J(1)w^{\ve}(z)|^ 2
]=0 ,
\end{equation}
which implies that the convergence of the finite-dimensional
distributions of $v^{\ve}$ can be reduced to
those of $w^{\ve}$. Hence, without loss of generality and from the
point of view
of the analysis we can assume that $\T=\mathrm{Id}$ and work with
\[
\nu^{\ve}(z)=J(1)\ve^{2h(z)-2}m\biggl({z\over\ve^2}, h(z)\biggr)=\ve
^{2h(z)-2}m\biggl({z\over\ve^2}, h(z)\biggr).
\]
Following \cite{CM}, the finite-dimensional distributions of the antiderivative
of $\nu^{\ve}$ converges to those of $S_h$, and thus ${\mathrm{A}}_1$ is
satisfied.
Now we check ${\mathrm{A}}_2(2)$.
We let $\delta>0$ and thanks to the asymptotic assumption on $m$,
there exists $z_{\delta}$ such that
for every~$z_1$, $z_2$ and $\ve$ satisfying $|z_1-z_2|>\ve^{2}z_{\delta
}$ we have
\begin{eqnarray*}
&&\sup_{(H_1, H_2)}\biggl| \biggl|{z_1-z_2\over\ve^2}\biggr|^{2-H_1-H_2}
\e\biggl[m\biggl({z_1\over\ve^{2}},H_1\biggr)
m\biggl({z_2\over\ve^{2}},H_2\biggr)\biggr]\\
&&\hspace*{133pt}{}-\mathcal{R}\biggl( {z_1\over\ve^{2}},{z_2\over\ve^{2}};H_1, H_2
\biggr) \biggr|<\delta.
\end{eqnarray*}
Then, noting that $\mathcal{R}( {z_1/ \ve^{2}},{z_2/ \ve^{2}},H_1, H_2) =
\mathcal{R}( z_1,z_2,H_1, H_2) $ and substituting $(H_1, H_2)$ by
$(h(z_1), h(z_2))$ we get
\begin{eqnarray*}
&&\biggl| \biggl|{z_1-z_2\over\ve^2}\biggr|^{2-h(z_1)-h(z_2)}\e\biggl[m
\biggl({z_1\over\ve^{2}},h(z_1)\biggr)
m\biggl({z_2\over\ve^{2}},h(z_2)\biggr)\biggr]\\
&&\hspace*{130pt}{}-\mathcal{R}( {z_1},{z_2};h(z_1),h(z_2)) \biggr|<\delta.
\end{eqnarray*}
Letting $\mathcal{R}^*( {z_1},{z_2}):=\mathcal{R}(
{z_1},{z_2};h(z_1),h(z_2))$ and noticing that
$\sup(1/\mathcal{R}^*)<\infty$ (because $\inf\mathcal{R}^*>0$) we obtain
\begin{eqnarray*}
&&\bigl| \e[\nu^{\ve}(z_1)\nu^{\ve}(z_2)] -
\mathcal{R}^*( {z_1},{z_2})|z_1-z_2|^{h(z_1)+h(z_2)-2}
\bigr|  \\
&&\quad{}<\delta\mathcal{R}^*( {z_1},{z_2})|z_1-z_2
|^{h(z_1)+h(z_2)-2}\sup
(1/\mathcal{R}^*),
\end{eqnarray*}
which proves ${\mathrm{A}_2}(2)$.
It remains to check ${\mathrm{A}_3}(2)$. Let $\rho>0$. Because of the
boundedness assumption
on $m$, there exists a
constant $C_1(\rho)>0$ so that
for every $z_1$, $z_2$ and $\ve$ satisfying $|z_1-z_2|/\ve^{2}< \rho$,
we have
\[
\biggl|\e\biggl[m\biggl({z_1\over\ve^{2}},H_1\biggr)
m\biggl({z_2\over\ve^{2}},H_2\biggr)\biggr]\biggr|\leq C_1(\rho).
\]
Thus,
\begin{eqnarray*}
| \e[\nu^{\ve}(z_1)\nu^{\ve}(z_2)]| & \leq& C_1(\rho)\ve
^{2h(z_1)+2h(z_2)-4}\\
& = & C_1(\rho)|z_1-z_2|^{h(z_1)+h(z_2)-2}
\biggl|{z_1-z_2\over\ve^2}
\biggr|^{2-h(z_1)-h(z_2)}\\
& \leq& C_1(\rho)|z_1-z_2|^{h(z_1)+h(z_2)-2} \rho
^{2-h(z_1)-h(z_2)}\\
& \leq& C_2(\rho)|z_1-z_2|^{h(z_1)+h(z_2)-2},
\end{eqnarray*}
where $C_2(\rho)$ can be chosen such that $C_1(\rho) \rho
^{2-h(z_1)-h(z_2)}\leq C_2(\rho)$.
So ${\mathrm{A}_3}(2)$ is satisfied and the proof can be concluded
by applying Theorem
\ref{main}.
\end{pf*}
We finish this subsection by applying Theorem \ref{main3} to an example
that was mentioned in \cite{CM}.
Let us consider $W_H$ defined as in (\ref{WH}).
We let
%
%e32 ###
\begin{equation}
\label{MH}
m(z, H)=W_H(z+1)-W_H(z).
\end{equation}
We compute the covariance between $m(z_1, H_1)$ and $m(z_2,H_2)$ for
every $z_1$, $z_2$, $H_1$ and $H_2$
%
%e33 ###
\begin{eqnarray}
\label{unif1}
&&\e[m(z_1, H_1)m(z_2,H_2)] \nonumber\\
 &&\quad= {1\over2}
{C({(H_1+H_2)}/{2} )^2\over C(H_1)C(H_2)}|z_1-z_2|^{H_1+H_2}
\\
&&\qquad{} \times\biggl( \biggl|1+{1\over z_1-z_2}\biggr|^{H_1+H_2}
+\biggl|1-{1\over z_1-z_2}\biggr|^{H_1+H_2}-2 \biggr).\nonumber
\end{eqnarray}
By Taylor's formula we get that the asymptotic covariance
$\R$ of $\{ m(z,H)\}_{z,H}$ can be written as
%
%e34 ###
\begin{equation}
\label{RH1H2}
\R(H_1,H_2)={1\over2}(H_1+H_2)(H_1+H_2-1){C( {(H_1+H_2)}/{2}
)^2\over C(H_1)C(H_2)} .
\end{equation}
Then applying
Theorem \ref{main3} we get that
$\{a^{\ve}(Z,s)\}_{s}$ converges
in distribution to $\widetilde{a}(Z,s) = f(s-{1\over2}S_h(Z)
)   $
where
%
%e35 ###
\begin{equation}
\label{rh}
{S}_h(Z) = J(1)\int_{-\infty}^{\infty} \biggl( \int_0^Z { -ix e^{-iux}
\over C(h(u))|x|^{h(u)+1/2} }du\biggr) \hat{B}(dx).
\end{equation}
As mentioned in Section 6.1 of \cite{CM}, we also can observe
that if we assume that $h$ is differentiable then we can write
${S}_h(Z)$ as
%
%e36 ###
\begin{eqnarray}
\label{rhbis}
& & {S}_h(Z) = J(1)\int_{-\infty}^{\infty} \hat{B}(dx) \biggl\{
{(e^{-iZx}-1)\over C(h(Z))
|x|^{h(Z)+1/2}}\nonumber\\
& &\hspace*{108pt}\quad {}  -\int_0^Z{(e^{-iu x}-1)\over|x|^{h(u
)+1/2}}\biggl(\frac{ \log|x|}{C(h(u))}\\
&&\hspace*{187pt}\qquad{}-\frac{C'(h(u))}{C(h(u))^2} \biggr) h'(u)  \, du \biggr\}
,\nonumber
\end{eqnarray}
which means that ${S}_h(Z)$ is the sum of a multifractional Brownian
motion as in~(\ref{Wh})
and of a regular process.

%s6 ###
\section{A non-Gaussian and multifractal medium}
\label{nonGm}
In this section we study the case of a medium that generalizes the
media discussed above.
We define $\{ m(z,H)\}_{z,H}$ for every $z\geq0$ by
%
%e37 ###
\begin{equation}
m(z,H)={1\over C(H)}\int_{\re}\exp(-izx){ \psi(x)
|x|^{1/2-H}}\hat{B}(dx),
\end{equation}
where $H\in(1/2, 1) $, $C(H)$ is a renormalization constant, $\psi$ is a
complex-valued symmetric function
and
$\hat{B}(dx)$ is the Fourier transform of a real Gaussian measure.
We assume that $\psi$ is continuous, $\psi(0)=1$ and satisfies $|\psi
(x)|=\mathcal{O}_{|x|\to\infty}(|x|^{-1})$.
Notice that the family of processes defined by (\ref{MH}) in terms of
fractional Brownian motion $\{ W_H(z) \}_{z,H}$ is an example of such a
process.

Thus, $\{ m(z,H)\}_{z,H}$ is a centered Gaussian field and its
covariance can be written as
%
%e38 ###
\begin{equation}
\label{unifx}
\e[m(z_1,H_1)m(z_2,H_2)]=\int_{\re}{\exp(i(z_2-z_1)x) |\psi(x)|^2
\over C(H_1)C(H_2)|x|^{H_1+H_2-1}}\,dx.
\end{equation}
Now we consider a function $h$ that takes its values in $[h_-,
h_+]\subset(1/2, 1)$ and a truncation
function $\T$ with Hermite index $K\in\mathbb{N}^*$. We define $\nu
^{\ve}$ as
\[
\nu^{\ve}(z)=\ve^{\kappa(z)-\tau}\nu\biggl( {z\over\ve^2 }, z \biggr),
\]
where
\[
\nu( z_1, z_2 )=\T( m ( z_1, \widetilde{h}_K(z_2)
) )
\]
with
\[
\widetilde{h}_K(z)={h(z)-1\over K}+1 .
\]
We can then show that $\nu^{\ve}$ satisfies assumptions ${\mathrm{A}}_2(2)$ and ${\mathrm{A}}_3(2)$.
In particular, we have
%
%e39 ###
\begin{equation}
\e[\nu^{\ve}(z_1)\nu^{\ve}(z_2)]\sim{J(K)^2\over{K!}}\R(\widetilde
{h}_K(z_1), \widetilde{h}_K(z_2))|z_1-z_2|^{h(z_1) + h(z_2)-2}
\end{equation}
when $|z_1-z_2|/\ve^2$ goes to $\infty$ assuming that
$\kappa(z)-\tau=2h(z)-2$, and $\R$ is defined as in (\ref{RH1H2}).
Therefore, because Theorem \ref{main} says that, under long-range
assumptions, the asymptotic behavior of
$a^{\ve}(Z, s)$ is essentially given by the limit of $v^{\ve}(z)$, we
can conclude by the following result.
\begin{theorem}
\label{main4}
As $\ve$ goes to 0, $\{a^{\ve}(Z,s)\}_{s}$ converges
in distribution in the space of continuous functions endowed with the
uniform topology to the random process $\{\widetilde{a}(Z,s)\}_{s}$
that can be written as
%
%e40 ###
\begin{equation}
\widetilde{a}(Z,s) = f\bigl(s-\tfrac{1}{2}S_h^K(Z) \bigr)   ,
\end{equation}
where
${S}_h^K$ is a centered process given for every $z$ by
%
%e41 ###
\begin{equation}
S_h^K(z)=
\int_{ \re^K} \mathcal{G}_{h,K}(z,x_1,\ldots, x_K ) \prod
_{k=1}^K\hat{B}(dx_k),
\end{equation}
where
\[
\mathcal{G}_{h,K}(z,x_1,\ldots, x_K )=\int_0^z \frac{J(K)e^{-iu\sum
_{k=1}^Kx_k}}{K! C(\widetilde{h}_K(u))^K} \prod_{k=1}^K
{-ix_k\over|x_k|^{\widetilde{h}_K(u)+1/2}}\,  du.
\]
\end{theorem}

Notice that the process $S_h^K$ is equal (in distribution) to $W_H^K$
of Section 4 if $h$ is a constant equal to $H$, and is equal
to $S_h$ of Section 5 if $K=1$. Because of these facts, $S_h^K$ is in
general non-Gaussian and multifractional. This shows that under general
long-range assumptions the asymptotic time-shift is neither Gaussian,
nor homogeneous. This is in dramatic contrast to the short-range case
where the time shift is a Brownian motion, which is homogeneous and Gaussian.

\begin{pf*}{Proof of Theorem \protect\ref{main4}}
We let
\[
v^{\ve}(z)  =  \int_0^z\nu^{\ve}(u) \, du = \int_0^z du\,
\ve^{2h(u)-2}\T\biggl( m \biggl( {u\over\ve^2}, \widetilde{h}_K(u)
\biggr) \biggr)
\]
and
\[
w_K^{\ve}(z)  =  \int_0^z du \,  \ve^{2h(u)-2} P_K
\biggl( m \biggl( {u\over\ve^2},\widetilde{h}_K(u) \biggr) \biggr) .
\]
Using the same arguments as for the beginning of the proof of Theorem
\ref{main3}, and
the fact that the Hermite index of $\T$ is $K$, we get
%
%e42 ###
\begin{equation}
\lim_{\ve\to0}\e\biggl[\biggl|v^{\ve}(z)-{J(K)\over K!}w_K^{\ve}(z)
\biggr|^ 2\biggr]=0.
\end{equation}
Then using the formula (see \cite{ito}, for instance)
\[
P_K\biggl(\int_{\re}\phi(x)\hat{B}(dx)\biggr)
=\int_{\re^K}\prod_{k=1}^K\phi(x_k)\hat{B}(dx_k)
\]
for every $\phi\in L^2(\re)$ we get
\begin{eqnarray*}
w_K^{\ve}(z) & = & \int_0^z du \,  \frac{\ve^{2h(u)-2}}{C(h(u))^K}
\int_{ \re^K}e^{-iu\sum_{j=1}^Kx_j/\ve^2}\prod_{k=1}^K{\psi(x_k)\over
|x_k|^{\widetilde{h}_K(u)-1/2}}\hat{B}(dx_k) \\ \nonumber
& = & \int_{ \re^K} \int_0^z du\,   \prod_{k=1}^K{\psi(x_k)\over
|x_k|^{\widetilde{h}_K(u)-1/2}}\hat{B}(dx_k) \frac{\ve^{2h(u)-2}}{C(h(u))^K}
e^{-iu\sum_{j=1}^Kx_j/\ve^2} .
\end{eqnarray*}
Then we make the substitution $x_k\to\ve^2x_k$ for every $k$
\begin{eqnarray*}
w_K^{\ve}(z) & = & \int_{ \re^K} \int_0^z du \,  \prod_{k=1}^K{\psi(\ve
^2x_k)\over
|\ve^2x_k|^{\widetilde{h}_K(u)-1/2}}\hat{B}(\ve^2\,dx_k)\frac{\ve
^{2h(u)-2}}{C(h(u))^K}
e^{-iu\sum_{j=1}^Kx_j}\\ \nonumber
& = & \ve^{-K}\int_{ \re^K} \int_0^z du \,  \prod_{k=1}^K{\psi(\ve
^2x_k)\over|x_k|^{\widetilde{h}_K(u)-1/2}}\hat{B}(\ve^2\,dx_k) \frac
{1}{C(h(u))^K}
e^{-iu\sum_{j=1}^Kx_j}.
\end{eqnarray*}
We let
\begin{eqnarray*}
\widetilde{w}_K^{\ve}(z) & = & \int_{ \re^K} \int_0^z du\,   \prod
_{k=1}^K{\psi(\ve^2x_k)\over
|x_k|^{\widetilde{h}_K(u)-1/2}}\hat{B}(dx_k) \frac{1}{C(h(u))^K}
e^{-iu\sum_{j=1}^Kx_j}.
\end{eqnarray*}
The self-similarity of the Brownian motion gives that $\hat{B}(\ve
^2\,dx_k)$ is equal in
distribution to $\ve\hat{B}(dx_k)$, then we get that
\[
{w}_K^{\ve}  \stackrel{\mathrm{f.d.d.}}{=}  \widetilde{w}_K^{\ve},
\]
where $\stackrel{\mathrm{f.d.d.}}{=}$ means the equality of the finite-dimensional distributions.
Then, using the assumptions on $\psi$, we obtain the
convergence a.s. of the finite-dimensional margins of ${J(K)\over
K!}\widetilde{w}_K^{\ve}$ to those
of $S_h^K$, and thus the convergence of the finite-dimensional
distributions of $v^{\ve}$
to those of $S_h^K$, so ${\mathrm{A}}_1$ is satisfied. Now, as observed
at the beginning of this section, using (\ref{unifx}) and by the same
procedure as in the
proof of Theorem \ref{main3} we show that ${\mathrm{A}}_2(2)$ and ${\mathrm{A}}_3(2)$ hold.
We can then conclude by Theorem
\ref{main}.
\end{pf*}

%f1 ###
\begin{figure}[b]

\includegraphics{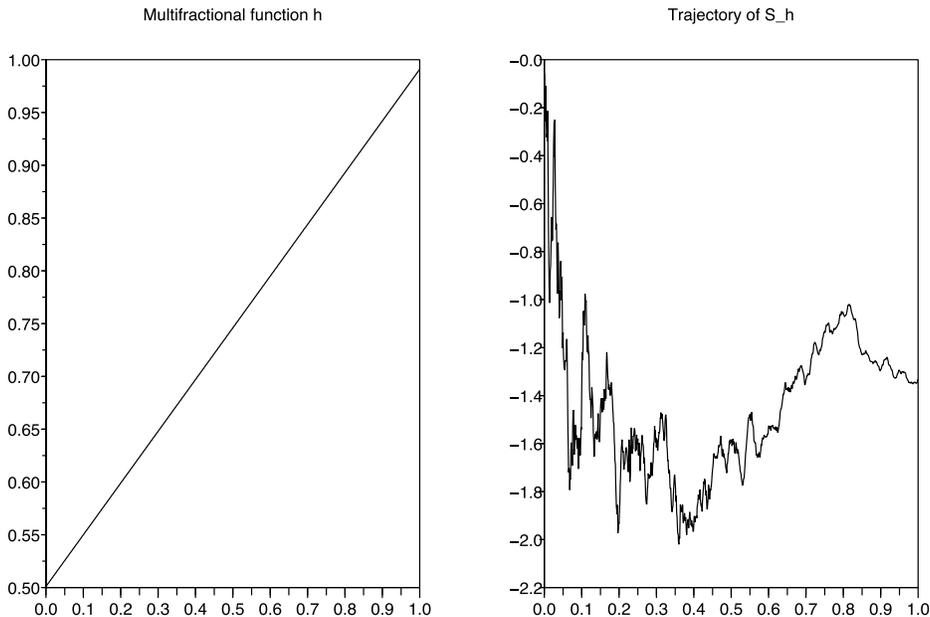}

\caption{Trajectory of $S_h$ with an increasing
multifractional function.}\label{courbe1}
\end{figure}

%s7 ###
\section{Numerical illustration}\label{secnum}
We illustrate our results with some numerical simulations.
In order to show the differences between the mixing and the
long-range cases, numerical simulations of the transmitted pulses
centered around the travel time are
presented in \cite{MS}. They are carried out
with a fractional white noise medium with Hurst index $H=0.5$ (corresponding
to the mixing case) and $H=0.6$ (corresponding to the long-range case).
These examples illustrate that the pulse shape is not
affected by the random fluctuation of the medium when $H=0.6$ and that
it is modified via a convolution with a Gaussian kernel when $H=0.5$.

\begin{figure}

\includegraphics{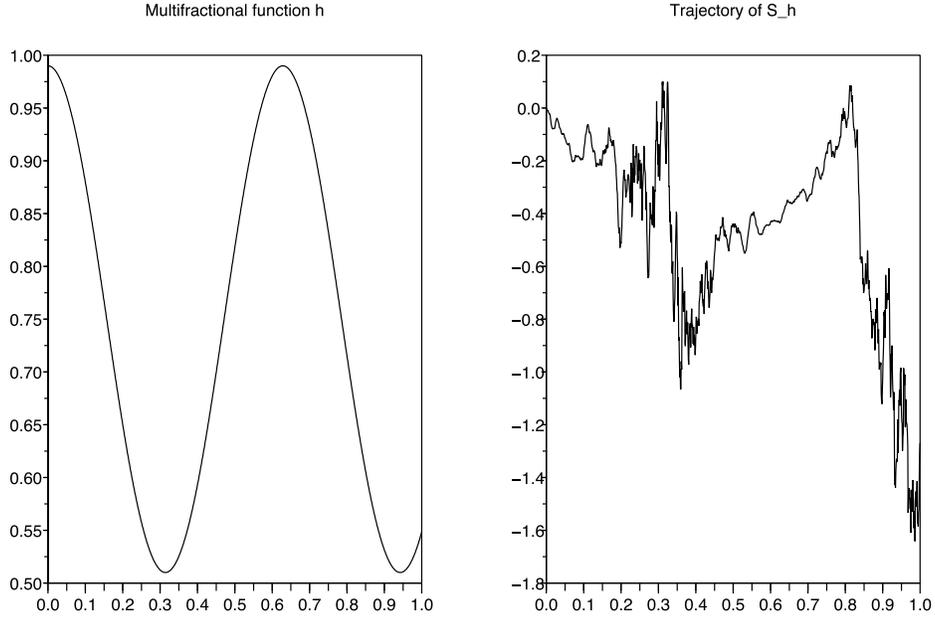}

\caption{Trajectory of $S_h$ with a periodic
multifractional function.}\label{courbe2}
\end{figure}

Here we aim to illustrate the differences between fractional and
multifractional cases. We present simulations of the asymptotic travel
times we obtain for
media with long-range correlation and different multifractional functions.
For the sake of simplicity we restrict ourself to the Gaussian case
presented in
Section~\ref{multifractal}, and we
let the propagation distance be one. For a fixed multifractional
function~$h$, the method we use to simulate the asymptotic travel time
$S_h$ is based on the method presented in \cite{T} (pages 370--371) and
the invariance principle proved in \cite{CM}. We first simulate the
fractional white noise $\{Y_j(H)\}_j$ of index $H\in(1/2,1)$ as in
equation (7.11.1) of
\cite{T} (page 371). Then, using Theorem 2 of \cite{CM} we can use
$\sum_{j=1}^{[Nt]}N^{-h(j/N)}Y_j(h(j/N))$ to approximate $S_h(t)$. In
Figure \ref{courbe1} we show a trajectory of $S_h$ with an increasing
multifractional function. In Figure \ref{courbe2} we show a trajectory of $S_h$
with an periodic multifractional function. In both figures we can
observe that the regularity varies along the trajectory
according to the local Hurst index. Modeling of this kind may, for instance,
be relevant in the case the multiscale crust of the sedimentary earth or
in the context of the turbulent atmosphere. In both cases the field is
typically strongly
anisotropic with a roughness that depends on depth or height, respectively.

%f2 ###

%s8 ###
\section[Proof of Theorem 3.1]{Proof of Theorem \protect\ref{main}} \label{secmain}
We first give an outline of the proof. As recalled in Section \ref{sec1} the process
$\{a^{\ve}(Z,s)\}_{s}$ can be written in
terms of the propagator $P_{\omega}^{\ve}$, and thus the study of the
convergence
of $\{a^{\ve}(Z,s)\}_{s}$ can be
analyzed via asymptotic properties of $P_{\omega}^{\ve}$.
The propagator $P_{\omega}^{\ve}$ satisfies the equation
\[
{dP_{\omega}^{\ve}\over dz}(z) =
\H_{\omega}^{\ve}
\biggl({z\over\ve^{\tau}},{z}\biggr)P_{\omega}^{\ve}(z)   ,
\]
that we can write in the form
%
%e43 ###
\begin{equation}
\label{propagator}
{dP_{\omega}^{\ve}(z)} =
{i\omega\over2}\sum_{j=1}^3F_jP_{\omega}^{\ve}(z) \,dv_j^{\ve}(z)   ,
\end{equation}
where
\[
F_1=\pmatrix{
1 & 0 \cr
0 & -1},\qquad
F_2=\pmatrix{
0 & -1 \cr 1 & 0}
\]
and
\[
F_3=\pmatrix{
0 & i \cr i & 0},
\]
and $v_1^{\ve}$, $v_2^{\ve}$ and $v_3^{\ve}$ are three processes of
bounded variation that we can write as
\begin{eqnarray*}
v_1^{\ve}(z)&=&\int_0^z
\nu^{\ve}({z'})  \,dz'  ,\\
v_2^{\ve}(z)&=&\int_0^z
\nu^{\ve}({z'})
\cos\biggl(2\omega{z'\over\ve^{\tau}}\biggr)  \,dz'  ,\\
v_3^{\ve}(z)&=&\int_0^z
\nu^{\ve}({z'})\sin\biggl(2\omega{z'\over\ve^{\tau}}
\biggr)
\,dz'   .
\end{eqnarray*}
Thanks to T. Lyons's rough paths theory for which we recall some tools in
the \hyperref[app]{Appendix} we shall see that
the convergence of $P_{\omega}^{\ve}$ can be reduced for a convenient
topology to the convergence of the
process $\mathbf{v}^{\ve}$ defined as
\[
\mathbf{v}^{\ve}:=(v_1^{\ve},v_2^{\ve},v_3^{\ve})  .
\]
Hence, we first prove the convergence of $\mathbf{v}^{\ve}$, then by
Theorem \ref{continuity}
(see the \hyperref[app]{Appendix}) we deduce the convergence of $P_{\omega}^{\ve}$ in
Section \ref{sec8.1} and thanks to (\ref{representation-a-fourier1}),
we finally conclude by the convergence of $\{a^{\ve}(Z,s)\}_{s}$ in
Section \ref{sec8.2}.

%s8.1 ###
\subsection{Convergence of the propagator}\label{sec8.1}
Using Theorem \ref{continuity} and the expression~(\ref{propagator}),
the asymptotic study
of the propagator is reduced to finding the limit in a rough path space of
$\mathbf{v}^{\ve}:=(v_1^{\ve},v_2^{\ve},v_3^{\ve})  $.
This is the aim of the following lemma.
\begin{lemma}
\label{convergence-v}
There exists $\gamma_*\in(0,1)$ such that for every $p>2/(2-\gamma
_*)$, as $\ve$ goes to 0, the increments of $\mathbf{v}^{\ve}$
converge in $\Omega_p$ to those of $\mathbf{V}$ which can be written as
\[
\mathbf{V}=(V, 0,0).
\]
\end{lemma}

The proof of Lemma \ref{convergence-v} is based on establishing several
technical lemmas that we do next. We let $r_{\nu^{\ve}}(x,y)=\e[\nu
^{\ve}(x)\nu^{\ve}(y)]$.
\begin{lemma}
\label{est1}
There exist $C$ and $\gamma_*$ so that
\[
|r_{\nu^{\ve}}(x,y)|\leq C|x-y|^{-\gamma_*}
\]
for every $x$ and $y$.
\end{lemma}
\begin{pf}
The assumptions of Theorem \ref{main} imply
that for every $\delta> 0$ there exists $z_{\delta}>0$ such that for
$|x-y|>\ve^{\lambda}z_{\delta}$ we have
\[
(1-\delta)R(x,y)|x-y|^{-\gamma(x,y)}\leq r_{\nu^{\ve}}(x,y)\leq
(1+\delta)R(x,y)|x-y|^{-\gamma(x,y)}.
\]
Hence, taking $\delta=1$ we get that for $|x-y|>\ve^{\lambda}z_{1}$ we have
\[
0\leq r_{\nu^{\ve}}(x,y)\leq C|x-y|^{-\gamma_+}.
\]
Moreover, thanks to the assumptions of Theorem \ref{main}, we know that
there exist $C_{z_1}$ and $\gamma_{z_1}$ so that
for $|x-y|\leq\ve^{\lambda}z_{1}$ we have
\[
0\leq|r_{\nu^{\ve}}(x,y)|\leq C_{z_1}|x-y|^{-\gamma_{z_1}}.
\]
By choosing $\gamma_*:=\max(\gamma_+, \gamma_{z_1})$ we get that there
exists $\gamma_*$ so that
\[
|\e[\nu^{\ve}(x)\nu^{\ve}(y)]|\leq C|x-y|^{-\gamma_*}
\]
for every $x$ and $y$.
\end{pf}
\begin{lemma}
\label{l3}
For every $z\in[0,Z]$, as $\ve$ goes to 0 the sequences $v_2^{\ve
}(z)$ and $v_3^{\ve}(z)$ converge to 0.
\end{lemma}
\begin{pf}
Without loss of generality we present the proof only for $v_2^{\ve}(z)$
and with
$2\omega=1$. We have
\begin{eqnarray*}
\e[v^{\ve}_2(z)^2] & = &
\int_0^zdx\int_0^zdy\,
\cos\biggl({x\over\ve^{\tau}}\biggr)\,\cos\biggl({y\over\ve^{\tau}}\biggr)
r_{\nu^{\ve}}({x}, {y})\\
& = & I_1^{\ve}(z)+I_2^{\ve}(z)   ,
\end{eqnarray*}
with
\begin{eqnarray*}
 I_1^{\ve}(z)&=&\int_0^zdx\int_0^zdy\,
\cos\biggl({x\over\ve^{\tau}}\biggr)\cos\biggl({y\over\ve^{\tau}}\biggr)
R({x}, {y})|{x-y}|^{-\gamma(x,y)}  ,\\
 I_2^{\ve}(z)&=&
\int_0^zdx\int_0^zdy\,
\cos\biggl({x\over\ve^{\tau}}\biggr)
\cos\biggl({y\over\ve^{\tau}}\biggr)\bigl(r_{\nu^{\ve}}({x,y})
-R(x,y) |{x-y}|^{-\gamma(x,y)}\bigr).
\end{eqnarray*}
Let $\delta>0$, and because of the assumptions of Theorem \ref{main},
we have that for $|x-y|>\ve^{\lambda}z_{\delta}$ (with $z_{\delta}$
sufficiently large)
$|r_{\nu^{\ve}}(x,y)-R(x,y)|x-y|^{-\gamma(x,y)}|\leq\delta
R(x,y)|x-y|^{-\gamma(x,y)}$ for every $\ve$.
Combining this with Lemma \ref{est1} we obtain
\begin{eqnarray*}
|I_2^{\ve}(z)|& \leq& {\delta}\int_0^zdx
\int_0^zdy\,R(x,y)
|{x-y}|^{-\gamma(x,y)}\\
&&{}+C_{\delta}\int_0^zdx\int_0^zdy\,
|x-y|^{-\gamma^* }1_{|{x-y}|\leq\ve^{\lambda}z_{\delta}}
\end{eqnarray*}
so that
\[
\limsup_{\ve\to0} |I_2^{\ve}(z)|\leq\delta\int_0^zdx\int
_0^zdy\,|{x-y}|^{-\gamma(x,y)}
.
\]
The inequality above is valid for every $\delta>0$, and we conclude
\[
\lim_{\ve\to0} I_2^{\ve}(z)=0  .
\]
We can deal with $ I_1^{\ve}(z)$ using a Riemann-type result. Indeed,
the function
$\tilde{R}\dvtx (x,y)\mapsto R(x,y)|x-y|^{-\gamma(x,y)}$ is integrable on
$\Delta_z=[0,z]^2$,
so we can approximate it by a sequence of constant by step functions
$(R_N)_N$ such that
\[
\lim_{N\to\infty}\int_0^zdx\int_0^zdy\,|\tilde{R}(x,y)-R_N(x,y)|=0.
\]
Moreover, we can write
\begin{eqnarray*}
|I_1^{\ve}(z)|&\leq&\biggl|\int_0^zdx\int_0^zdy\,
\cos\biggl({x\over\ve^{\tau}}\biggr)
\cos\biggl({y\over\ve^{\tau}}\biggr)R_N(x,y)\biggr|\\
&&{}+\int_0^zdx\int
_0^zdy\,|\tilde{R}(x,y)-R_N(x,y)|
\end{eqnarray*}
for every $\ve$ and $N$.
We easily see that
\[
\lim_{\ve\to0}\int_0^zdx\int_0^zdy\,
\cos\biggl({x\over\ve^{\tau}}\biggr)
\cos\biggl({y\over\ve^{\tau}}\biggr)R_N(x,y)=0
\]
so that
\[
\limsup_{\ve\to0}
|I_1^{\ve}(z)|\leq\int_0^zdx\int_0^zdy\,|\tilde{R}(x,y)-R_N(x,y)|
\]
for every $N$.
This finally shows
\[
\lim_{\ve\to0} I_1^{\ve}(z)=0
\]
and then
\[
\lim_{\ve\to0} \e[v^{\ve}_2(z)^2]=0   ,
\]
which completes the proof.
\end{pf}

Now we deal with a technical lemma regarding the increments of $\mathbf
{v}^{\ve}$.
\begin{lemma}
\label{l5}
There exist $C>0$ and $\gamma_*\in(0,1)$ such that for every $z$,
$\zeta$ and $\ve>0$ we have
\[
\e[\| \mathbf{v}^{\ve}(z)-\mathbf{v}^{\ve}({\zeta}) \|^2]\leq C|z-\zeta
|^{2-\gamma^*}.
\]
\end{lemma}
\begin{pf}
Because of Lemma \ref{est1} there exists $\gamma_*$ so that
$|\e[\nu^{\ve}(x)\nu^{\ve}(y)]|\leq C|x-y|^{-\gamma_*}$ for every $x$
and $y$.
Then, for every $j=1,2,3$, we have (taking $z>{\zeta}$)
\begin{eqnarray*}
\e[| v_j^{\ve}(z)-v_j^{\ve}({\zeta}) |^2] & \leq&
\int_{\zeta}^zdx\int_{\zeta}^zdy \,|\e[ \nu^{\ve}({x}
)\nu^{\ve}({y}) ]|\\
& \leq& {C}\int_{\zeta}^zdx\int_{\zeta}^zdy\,|{x-y}|^{-\gamma
^*} \\
& \leq& \frac{2C'}{(1-\gamma^*)(2-\gamma^*)}|z-\zeta
|^{2-\gamma^*},
\end{eqnarray*}
which completes the proof.
\end{pf}

In the sequel we shall use the
notation $H_*:=(2-\gamma_*)/2$.
Using the above lemmas we next deduce the following lemma which deals
with identification of the limit.
\begin{lemma}
\label{l6}
The process $\mathbf{V}$ defined in Lemma \ref{convergence-v} is a.s.
continuous (up to a modification).
Moreover, as $\ve$ goes to 0, $\mathbf{v}^{\ve}$ converges to $\mathbf{V}$
in the space of continuous functions endowed with the uniform norm.
\end{lemma}
\begin{pf}
Assumptions and Lemma \ref{l3} give the convergence of finite-\break dimensional distributions
of $\mathbf{v}^{\ve}$ to those of $\mathbf{V}$.
Using then the Kolmogorov criterion~\cite{B},
Lemma \ref{l5}, and the fact that $2H_*>1$ we get
the tightness of $(\mathbf{v}^{\ve})_{\ve}$ in
the space of continuous functions
endowed with the uniform norm which establishes the proof.
\end{pf}%\hspace*{298pt}

Thanks to Lemma \ref{l6} we conclude with the proof of Lemma \ref
{convergence-v} by establishing the tightness in a rough paths sense.
\begin{lemma}\label{l7}
For every $p>1/H_*$, the sequence $(\mathbf{v}^{\ve})_{\ve}$ is tight
in $\Omega_p$ and the process $\mathbf{V}$ is a.s. of finite $p$-variation.
\end{lemma}
\begin{pf*}{Proof of Lemmas \protect\ref{convergence-v} and \protect\ref{l7}}
Let $q\in(1/H_*,p)$. In view of Lemmas \ref{lejaybis} and \ref{l6} it
is enough to prove
%
%e44 ###
\begin{equation}
\label{l7p}
\lim_{A\to+\infty}\sup_{\ve
>0}\mathbb{P}[V_q(\mathbf{v}^{\ve})>A]=0  .
\end{equation}
Using Chebyshev's inequality, the fact that $q<2$, Lemma \ref{ledoux},
the H\"older inequality and Lemma \ref{l6}
we find
\begin{eqnarray*}
\mathbb{P}[V_q(\mathbf{v}^{\ve})>A] & \leq& {1\over A^q}\e[V_q(\mathbf
{v}^{\ve})^q]\\
& \leq& {C\over A^q}\sum_{n=1}^{+\infty}n^{C}\sum_{k=1}^{2^n}\e[\|
\mathbf{v}^{\ve}(z_{k}^n)- \mathbf{v}^{\ve}(z_{k-1}^n)\|^q]\\
& \leq& {C\over A^q}\sum_{n=1}^{+\infty}n^{C}\sum_{k=1}^{2^n}\e[\|
\mathbf{v}^{\ve}(z_{k}^n)- \mathbf{v}^{\ve}(z_{k-1}^n)\|^2]^{q/2}\\
& \leq& {C'\over A^q}\sum_{n=1}^{+\infty}n^{C}\sum_{k=1}^{2^n}{
\biggl({1\over2^n}\biggr)}^{qH_*}\\
& \leq& {C'\over A^q}\sum_{n=1}^{+\infty}n^{C}{\biggl({1\over
2^n}\biggr)}^{qH_*-1}  ,
\end{eqnarray*}
and since $qH_*>1$ we deduce (\ref{l7p}).
\end{pf*}

Finally, we can now derive the following lemma which deals with the
convergence of the propagator.
\begin{lemma}
\label{conv-prop}
Let $\{ \omega_1, \ldots, \omega_n \}$ to be a collection of
frequencies. Then, as $\ve$ goes to 0, the propagator vector $(P_{\omega
_1}^{\ve}, \ldots, P_{\omega_n}^{\ve})$ converges in distribution in
the space of continuous
functions to $(P_{\omega_1}, \ldots, P_{\omega_n})$ which is the
asymptotic propagator
$P_{\omega}$ that we can write as
\[
P_{\omega}(z)
=
\pmatrix{
\exp\biggl(\dfrac{i\omega}{2}{V}(z)\biggr) & 0\vspace*{2pt}\cr
0 &  {\exp\biggl(-\dfrac{i\omega}{2}{V}(z)\biggr)}}
.
\]
\end{lemma}
\begin{pf}
By combining Theorem \ref{continuity}, (\ref{propagator}) and Lemma \ref
{convergence-v} we get
that, as~$\ve$ goes to 0, $P^{\ve}_{\omega}$ converges in distribution
in the space of continuous
functions (endowed with the topology of the uniform convergence)
to the solution $P_{\omega}$ of the following system of equations:
\begin{eqnarray*}
dP_{\omega}(z)={i\omega\over2}
\pmatrix{1 & 0 \cr 0 & -1}P_{\omega}(z)\,   d{V}(z)  .
\end{eqnarray*}
This concludes the proof.
\end{pf}

We remark that the situation here contrasts with the short-range case.
Indeed, the asymptotic propagator is driven by one process in the
long-range case whereas it is driven by three processes in the short-range case.

%s8.2 ###
\subsection{Conclusion of the proof}\label{sec8.2}
The remaining part of the proof of Theorem \ref{main} follows the lines of
\cite{CF,book1};
however, we present it here for completeness.
Recall that thanks to the formula (\ref{representation-a-fourier1}) we
can write $a^{\ve}(Z,s)$ in a Fourier-type formula using the
transmission coefficient
%
%e45 ###
\begin{eqnarray}
\label{representation-a-fourier-bis}
a^{\ve}(Z,s) ={1\over2\pi}\int e^{-is\omega}T_{\omega}^{\ve}(Z)\widehat
{f}(\omega)\,  d\omega
  ,
\end{eqnarray}
with the transmission coefficient being a functional of the propagator
$P^{\ve}_{\omega}$.
We shall use Lemma \ref{conv-prop} to deduce the convergence of the
transmitted wave.

Let $n\in\mathbb{N}$, $s_1\leq\cdots \leq s_n\in[0,\infty)$.
We can write:
\begin{eqnarray*}
&& \e[a^{\ve}(Z,s_1)\cdots a^{\ve}(Z,s_n) ]\\
&&\quad = \e\Biggl[ {1\over
(2\pi)^{n}}\prod_{j=1}^{n}\int e^{-is_j\omega}T_{\omega}^{\ve
}(Z)\hat{f}(\omega) \, d\omega\Biggr] \\
&&\quad = {1\over(2\pi)^{n}}
\int\cdots\int e^{-i\sum_{j=1}^{n} s_j\omega_j}\hat{f}(\omega
_1)\cdots\hat{f}(\omega_n)\\
&&\hspace*{64pt}\qquad{}\times\e[T_{\omega_1}^{\ve}(Z)\cdots T_{\omega_n}^{\ve}(Z) ]  \,d\omega
_1\,\cdots\, d\omega_n  .
\end{eqnarray*}
Thanks to Lemma \ref{conv-prop} we have that as $\ve\to0$
\[
\e[T_{\omega_1}^{\ve}(Z)\cdots T_{\omega_n}^{\ve}(Z) ]\longrightarrow \e\Biggl[ \exp
\Biggl({i{V}(Z)\over2}\sum_{j=1}^{n}\omega_j\Biggr)\Biggr]   ,
\]
and then
\begin{eqnarray*}
&&\e[a^{\ve}(Z,s_1)\cdots a^{\ve}(Z,s_n) ]\\
&&\quad\to
{1\over(2\pi)^{n}}
\int\cdots\int e^{-i\sum_{j=1}^{n} s_j\omega_j}\hat{f}(\omega
_1)\cdots\hat{f}(\omega_n)
\\
&&\hspace*{67pt}\qquad{} \times
\e\Biggl[ \exp\Biggl({i{V}(Z)\over2}\sum_{j=1}^{n}\omega_j\Biggr)\Biggr]
\,d\omega_1\,\cdots \,d\omega_n
\\
&&\quad =\e\Biggl[ {1\over(2\pi)^{n}}\prod_{j=1}^{n}\int
e^{-i(s_j-{{V}(Z)/2})\omega}\hat{f}(\omega) \, d\omega\Biggr] \\
&&\quad  =\e\Biggl[ \prod_{j=1}^{n}f\bigl(s_j-{{V}(Z)/2}\bigr)\Biggr]   .
\end{eqnarray*}
The tightness proof is similar to the proof of Lemma 3.2 in \cite{CF}
and the convergence of $a^{\ve}(Z,s)$ follows.

\begin{appendix}\label{app}

%s9 ###
\section*{Appendix: Differential equations and rough paths}
In this appendix we fix $p\in[1, 2)$ and consider a closed interval $I=[0,Z]$.
We define the $p$-variation of a continuous function $w\dvtx I\to\re^n$ by
\[
V_p(w):=\Biggl( \sup_D\sum_{j=0}^{k-1}\| w(z_{j+1})-w(z_{j})\|^p
\Biggr)^{1/p}  ,
\]
where $\sup_D$ runs over all finite partition $\{ 0=z_0, \ldots, z_k=Z \}$
of $I$ and
where here and below $\| \cdot\|$ refers to the $L^2$ norm.
The space of all continuous functions of bounded
variation
(1-variation)
is endowed with the $p$-variation distance
\[
\|w\|_p=V_p(w)+\sup_{z\in
[0,Z]}|w(z)|   ,
\]
and is denoted by $\Omega^{\infty}_p$. The closure of this metric space
is called the space of all geometric rough paths and is denoted by
$\Omega_p$. One of the most important theorems of rough paths theory
is the following:
\begin{theorem}[(T. Lyons's Continuity Theorem)]\label
{continuity}
Let\footnote{Here $\mathcal{L}(\Rset, \Rset^d)$ [resp., $\mathcal
{L}(\Rset^n, \Rset^d)$]
denotes the space of all linear
maps from $\Rset$ (resp., $\Rset^n$) to $\Rset^d$.} $G\dvtx \Rset\times\Rset
^d\to\mathcal{L}(\Rset, \Rset^d)$
and $F\dvtx \Rset\times\Rset^d\to\mathcal{L}(\Rset^n,  \Rset^d)$ be two
smooth functions.
Let $y$ be the unique
solution of the differential equation
\[
dy(z)=G(z,y(z))  \,dz+F(z,y(z)) \,dw(z),\qquad     y(z=0)=y_0  ,
\]
where $w$ is a bounded variation function. Then It\^{o}'s map
$\mathcal{I}\dvtx {w}\mapsto{y}$ is continuous with respect to the
$p$-variation distance from $\Omega^{\infty}_p(\Rset^n)$ to
$\Omega^{\infty}_p(\Rset^d)$. Therefore there exists a unique extension
of this
map (that we still denote by $\mathcal{I}$) to the space $\Omega_p(\Rset^n)$
\end{theorem}

This theorem has been proved by T. Lyons and extensively studied and applied
(see \cite{CQ1,LLQ,Le1,L1}).

The proof of Theorem \ref{main} is based on
analysis of
the tightness in the space of
geometric rough paths. In the context of this
we need to compute the $p$-variation for
$p>1$.
To this effect we will need the following lemmas of which
the first can be found, for instance, in
\cite{Le1}, and the second in \cite{LLQ,Le1}.
\setcounter{lemma}{0}
\begin{lemma}
\label{lejaybis}
Let $q\in[1,2)$ and $(v^{\ve})_{\ve
>0}$ a
family of continuous random processes of finite $q$-variation whose
associated family of probability measures is tight in the space of
continuous functions on $I$ and satisfying
%
%e46 ###
\setcounter{equation}{0}
\begin{equation}
\label{lejay2bis}
\lim_{A\to+\infty}\sup_{\ve
>0}\mathbb{P}[V_q(v^{\ve})>A]=0  .
\end{equation}
Then the family of probability measures associated to $(v^{\ve})_{\ve>0}$
is tight in $\Omega_p$ for every $p>q$.
\end{lemma}
\begin{lemma}
\label{ledoux}
For every $n\in\mathbb{N}$ and every $k=0, 1,\ldots, 2^n$, we let
$z_k^n:=Z k/2^n$. Let $q\in[1,2)$ and $v$ be a function of finite
$q$-variation.
Then there exist two positive constants
$C_1$, $C_2$ which do not depend on $v$ such that
\[
V_q(v)^{q} \leq C_1\sum_{n=1}^{+\infty}n^{C_2}\sum_{k=1}^{2^n}\|
v(z_{k}^n)- v(z_{k-1}^n)\|^q  .
\]
\end{lemma}
\end{appendix}

%
% imsref loaded by akundreckaite, 2010-07-13 13:30:27

%

\printaddresses

\end{document}